\documentclass[10pt,a4paper,reqno]{amsart}
\usepackage{amsthm}
\usepackage{amsmath}
\usepackage{amsmath, amssymb}
\usepackage{type1cm}
\usepackage{cite}
\usepackage{bm}
\usepackage{cases}
\theoremstyle{definition}
\newtheorem{defi}{Definition}
\newtheorem{prop}{Proposition}
\newtheorem{lemma}{Lemma}
\newtheorem{theorem}{Theorem}
\newtheorem{corollary}{Corollary}
\newcommand{\eqdef}{\stackrel{\mathrm{def}}{=}}
\begin{document} 
\title{Generalized Bregman and Jensen divergences which include some f-divergences}
\author{Tomohiro Nishiyama}
\begin{abstract}
In this paper, we introduce new classes of divergences by extending the definitions of the Bregman divergence and the skew Jensen divergence. 
These new divergence classes (g-Bregman divergence and skew g-Jensen divergence) satisfy some properties similar to the Bregman or skew Jensen divergence.
We show these g-divergences include divergences which belong to a class of f-divergence (the Hellinger distance, the chi-square divergence and the alpha-divergence in addition to the Kullback-Leibler divergence).
Moreover, we derive an inequality between the g-Bregman divergence and the skew g-Jensen divergence and show this inequality is a generalization of Lin's inequality.
\\

\end{abstract}
\date{}
\maketitle
\bibliographystyle{plain}
\section{Introduction}
Divergences are functions measure the discrepancy between two points and play a key role in the field of machine learning, signal processing and so on.
Given a set $S$ and $p,q\in S$, a divergence is defined as a function $D:S\times S\rightarrow \mathbb{R}$ which satisfies the following properties.\\
1. $D(p,q)\geq 0 $ for all $p,q\in S$\\
2. $D(p,q)=0\iff p=q$\\

The Bregman divergence \cite{bregman1967relaxation} $B_{F}(p,q)\eqdef F(p)-F(q)-\langle p-q,\nabla F(q)\rangle$ and $f$-divergence \cite{csiszar1967information,ali1966general} $D_f(p\| q)\eqdef \sum_i q_i f(\frac{p_i}{q_i})$ are representative divergences defined by using a strictly convex function, where $F,f:S\rightarrow \mathbb{R}$ are strictly convex functions and $f(1)=0$.  
Besides that Nielsen have introduced the skew Jensen divergence $J_{F,\alpha}(p ,q)\eqdef (1-\alpha)F(p)+\alpha F(q)-F((1-\alpha)p+\alpha q)$ for a paramter $\alpha\in(0,1)$ and have shown the relation between the Bregman divergence and the skew Jensen divergence\cite{nielsen2011burbea, burbea1980convexity}.

The Bregman divergence and the f-divergece include well-known divergences.
For example, the Kullback-Leibler divergence(KL-divergence) \cite{kullback1997information} is a type of the Bregman divergence and the $f$-divergence.
On the other hand, the Hellinger distance, the $\chi$-square divergence and the $\alpha$-divergence \cite{cichocki2010families,amari2007methods} are a type of the $f$-divergence. 
The Jensen-Shannon divergence (JS-divergence)\cite{lin1991divergence} is a type of the skew Jensen divergence.

For probability distributions $p$ and $q$, each divergence is defined as follows.\\
\textbf{KL-divergence}\\
\begin{equation}
KL(p\| q)\eqdef \sum_i p_i\ln(\frac{p_i}{q_i})
\end{equation}
\textbf{Hellinger distance}\\
\begin{equation}
H^2(p, q)\eqdef \sum_i (\sqrt{p_i}-\sqrt{q_i})^2\\
\end{equation}
\textbf{Pearson $\chi$-square divergence}\\
\begin{equation}
D_{\chi^2,P}(p\| q)\eqdef \sum_i \frac{(p_i-q_i)^2}{q_i}\\
\end{equation}
\textbf{Neyman $\chi$-square divergence}\\
\begin{equation}
D_{\chi^2,N}(p\| q)\eqdef \sum_i \frac{(p_i-q_i)^2}{p_i}\\
\end{equation}
\textbf{$\alpha$-divergence}\\
\begin{equation}
D_{\alpha}(p\| q)\eqdef \frac{1}{\alpha(\alpha-1)}\sum_i (p_i^\alpha q_i^{1-\alpha} -1)\\
\end{equation}
\textbf{JS-divergence}\\
\begin{equation}
JS(p, q)\eqdef  \frac{1}{2}KL\bigl(p\| \frac{p+q}{2}\bigr)+\frac{1}{2}KL\bigl(q\| \frac{p+q}{2}\bigr)
\end{equation}

The main goal of this paper is to introduce new classes of divergences (''$g$-divergence'')  and study properties of the $g$-divergences.

First, we introduce the $g$-Bregman divergence, the symmetric $g$-Bregman divergence and the skew $g$-Jensen divergence by extending the definitions of the Bregman divergence and the skew Jensen divergence respectively. This idea is based on Zhang's papers \cite{zhang2013nonparametric,zhang2004divergence}.
For the $g$-Bregman divergence, we show some geometrical properties (for triangle, parallelogram, centroids) and derive an inequality between the symmetric $g$-Bregman divergence and the skew $g$-Jensen divergence.

Then, we show that the $g$-Bregman divergence includes the Hellinger distance, the Pearson and Neyman $\chi$-square divergence and the $\alpha$-divergence. Furthermore, we show that the skew $g$-Jensen divergence include the Hellinger distance and the $\alpha$-divergence. These divergences have not only the properties of $f$-divergence but also some properties similar to the Bregman or the skew Jensen divergence. 

Finally we derive many inequalities by using an inequality between the $g$-Bregman divergence and skew $g$-Jensen divergence.

\section{The $g$-Bregman divergence and the skew $g$-Jensen divergence}%
\subsection{Definition of the $g$-Bregman divergence and the skew $g$-Jensen divergence}
By using a injective function $g$, we define the $g$-Bregman divergence and the skew $g$-Jensen divergence.

\begin{defi}
Let $\Omega\subset \mathbb{R}^d$ be a convex set and let $p,q$ be points in $\Omega$. 
Let  $F(x)$ be a strictly convex function  $F:\Omega\rightarrow \mathbb{R}$.

The Bregman divergence and the symmetric Bregman divergence are defined as 
\begin{align}
B_{F}(p ,q)\eqdef F(p)-F(q)-\langle p-q, \nabla F(q)\rangle \\
\label{def_sym_Bregman}
B_{F,\mathrm{sym}}(p,q)\eqdef B_{F}(p,q)+B_{F}(q,p)=\langle p-q, \nabla F(p)- \nabla F(q) \rangle.
\end{align}
\end{defi}
The symbol $\langle , \rangle$ indicates an inner product.

For the parameter $\alpha\in\mathbb{R}\setminus\{0,1\}$, the scaled skew Jensen divergence \cite{nielsen2011burbea,nielsen2010family } is defined as 
\begin{equation}
sJ_{F,\alpha}(p ,q)\eqdef \frac{1}{\alpha(1-\alpha)}\biggl((1-\alpha)F(p)+\alpha F(q)-F((1-\alpha)p+\alpha q)\biggr).
\end{equation}

Nielsen has shown the relation between the Bregman divergence and the Jensen divergence as follows \cite{nielsen2011burbea}.
\begin{align}
\label{jensen_lim1}
B_F(q,p)&=\lim_{\alpha\to 0}sJ_{F,\alpha}(p,q)\\
\label{jensen_lim2}
B_F(p,q)&=\lim_{\alpha\to 1}sJ_{F,\alpha}(p,q)
\end{align}
\begin{equation}
\label{Jensen_Bregman_eq}
sJ_{F,\alpha}(p,q)=\frac{1}{\alpha(1-\alpha)}\biggl((1-\alpha)B_F(p,(1-\alpha)p+\alpha q)+\alpha B_F(q,(1-\alpha)p+\alpha q)\biggr)
\end{equation}

\begin{defi}(Definition of the $g$-convex set)
Let $S$ be a vector space over the real numbers and $g$ be a injective function $g:S\rightarrow S$.
If $S$ satisfies the following condition, we call the set $S$ ''$g$-convex set''.

For all $p$ and $q$ in $S$ and all $\alpha$ in the interval $(0, 1)$, the point $g^{-1}\bigl((1-\alpha)g(p) + \alpha g(q)\bigr)$ also belongs to $S$.
\end{defi}

\begin{defi} (Definition of the $g$-divergences)
Let $\Omega$ be a $g$-convex set and let $g:\Omega\rightarrow\Omega$ be a function which satisfies $g(p)=g(q)\iff p=q$.
Let $F:\Omega\rightarrow\mathbb{R}$ be a strictly convex function.
Let $\alpha$ be a parameter in $(0,1)$.

We define the $g$-Bregman divergence, the $g$-symmetric Bregman divergence and the (scaled) skew $g$-Jensen divergence as follows. 
\begin{align}
B_{F}^g(p,q)\eqdef B_{F}(g(p),g(q))\\
\label{def_sym_breg}
B_{F,\mathrm{sym}}^g(p,q)\eqdef B_{F,\mathrm{sym}}(g(p),g(q))\\
sJ_{F,\alpha}^g(p ,q)\eqdef sJ_{F,\alpha}(g(p) ,g(q))
\end{align}
\end{defi}
From the definition of function $g$, these functions satisfy divergence properties $D^g(p,q)=0\iff g(p)=g(q)\iff p=q$ and $D^g(p, q)\geq 0$, where $D^g$ denotes the $g$-divergences.
When $g(p)$ is equal to $p$ for all $p\in\Omega$, the $g$-divergences are consistent with original divergences.

In the following, we assume that $g$ is a function having an inverse function.

\subsection{Geometrical properties of the $g$-Bregman divergence}
In this subsection, we show the $g$-Bregman divergence and the symmetric $g$-Bregman divergence satisfy some geometrical properties as well as the Bregman and the symmetric Bregman divergence.

First, we show the $g$-Bregman divergence and the symmetric $g$-Bregman divergence satisfy the
linearity, the generalized law of cosines and the generalized parallelogram law.

Then, we show the point which minimize the weighted average of $g$-Bregman divergence ($g$-Bregman centroids) is equal to the quasi-arithmetic mean\cite{porcu2009quasi, amari2009information}.

This property can be used for clustering algorithms such as k-means algorithm\cite{lloyd1982least}.
\\
\textbf{Linearity}\\
For positive constant $c_1$ and $c_2$, 
\begin{equation}
B_{c_1F_1+c_2F_2}^g(p,q)=c_1B_{F_1}^g(p,q)+c_2B_{F_2}^g(p,q)
\end{equation}
holds.

\begin{prop}\textbf{(Generalized law of cosines)}\\
\label{cor_cosine}
Let $\Omega$ be a $g$-convex set.
For points $p,q,r\in \Omega$, the following equations hold.
\begin{align}
\label{triangular relation}
B_{F}^g(p,q)=B_{F}^g(p,r)+B_{F}^g(r,q)-\langle g(p)-g(r), \nabla F(g(q))- \nabla F(g(r))\rangle \\
B_{F,\mathrm{sym}}^g(p,q)=B_{F,\mathrm{sym}}^g(p,r)+B_{F,\mathrm{sym}}^g(r,q) \\ \nonumber
-\langle g(p)-g(r), \nabla F(g(q))- \nabla F(g(r))\rangle+\langle g(q)-g(r), \nabla F(g(p))- \nabla F(g(r))\rangle
\end{align}
\end{prop}
It is easily proved in the same way as the Bregman divergence by putting $p'=g(p)$, $q'=g(q)$ and $r'=g(r)$.

\begin{theorem}\textbf{(Generalized parallelogram law)}\\
\label{parallelogram}
Let $\Omega$ be a $g$-convex set.
When points $p,q,r,s\in \Omega$ satisfy $g(p)+g(r)=g(q)+g(s)$, the following equations holds.
\begin{equation}
B_{F,\mathrm{sym}}^g(p,q)+B_{F,\mathrm{sym}}^g(q,r)+B_{F,\mathrm{sym}}^g(r,s)+B_{F,\mathrm{sym}}^g(s,p)=B_{F,\mathrm{sym}}^g(p,r)+B_{F,\mathrm{sym}}^g(q,s)
\end{equation}
\end{theorem}

The left hand side is the sum of four sides of rectangle $pqrs$ and the right hand
side is the sum of diagonal lines of rectangle $pqrs$.\\
\noindent\textbf{Proof.}
It has been shown that the Bregman divergence satisfies the four point identity(see equation (2.3) in \cite{reich2010two}).
\begin{equation}
\langle\nabla F(r)-\nabla F(s), p-q\rangle=B_F(q,r)+B_F(p,s)-B_F(p,r)-B_F(q,s)
\end{equation}
By putting $p'=g(p)$, $q'=g(q)$, $r'=g(r)$ and $s'=g(s)$, we can prove the $g$-Bregman divergence satisfies the similar four point identity.
\begin{equation}
\langle\nabla F(g(r))-\nabla F(g(s)), g(p)-g(q)\rangle=B_F^g(q,r)+B_F^g(p,s)-B_F^g(p,r)-B_F^g(q,s)
\end{equation}
By combining the assumption and the definition of the symmetric $g$-Bregman divergence ((\ref{def_sym_Bregman}) and (\ref{def_sym_breg})), we have
\begin{equation}
\label{four_point_identity}
-B_{F,\mathrm{sym}}^g(r,s)=B_F^g(q,r)+B_F^g(p,s)-B_F^g(p,r)-B_F^g(q,s).
\end{equation}
Exchanging $p\leftrightarrow r$ and $q\leftrightarrow s$ in (\ref{four_point_identity}) and taking the sum with  (\ref{four_point_identity}), the result follows.
This theorem can be also proved in the same way as mentioned in the paper\cite{nishiyama2018divergence}.

\begin{defi} (Definition of the multivariate skew $g$-Jensen divergence)
\label{definition of multivariate}
Let $\Omega$ be a $g$-convex set and let $g:\Omega\rightarrow\Omega$ be a function which satisfies $g(p)=g(q)\iff p=q$.
Let $F:\Omega\rightarrow\mathbb{R}$ be a strictly convex function.
Let $p_\nu(\nu=1,2,3\cdots,N)$ are points in $\Omega$.
Let $\alpha_\nu\geq 0(\nu=1,2,3\cdots,N)$ be parameters which satisfy $\sum_{\nu=1}^N\alpha_\nu=1$.

We define the multivariate skew $g$-Jensen divergence as follows.
\begin{align}
J_{F,\boldsymbol{\alpha}}^g(p_1,p_2,\cdots, p_N)\eqdef\sum_{\nu=1}^N \alpha_\nu F(g(p_\nu))-F\biggl( \sum_{\nu=1}^N \alpha_\nu g(p_\nu)\biggr),
\end{align}
where $\boldsymbol{\alpha}$ denotes a vector $(\alpha_1, \alpha_2, \cdots, \alpha_N)$.
\end{defi}

\begin{theorem}\textbf{(The centroids of the $g$-Bregman divergence)}\\
\label{centroid_theorem}
Let $\Omega$ be a $g$-convex set.
Let $p_\nu(\nu=1,2,3\cdots,N)$ and $q$ are points in $\Omega$.
Let $\alpha_\nu\geq 0(\nu=1,2,3\cdots,N)$ be parameters which satisfy $\sum_{\nu=1}^N\alpha_\nu=1$.

Then, the following inequality holds.
\begin{align}
\sum_{\nu=1}^N \alpha_\nu B_F^g(p_\nu,q)\geq J_{F,\boldsymbol{\alpha}}^g(p_1,p_2,\cdots, p_N)
\end{align}
Equality holds if and only if
\begin{align}
q=g^{-1}\biggl(\sum_{\nu=1}^N \alpha_\nu g(p_\nu) \biggr)
\end{align}
Because the function $g$ is injective by the definition of the $g$-Bregman divergence, this value is the quasi-arithmetic mean. 
For example, 
\begin{subnumcases}
{g(p)=}
p & ($q$: the arithmetic mean) \\
\ln p & ($q$: the geometric mean) \\
p^{-1} & ($q$: the harmonic mean)\\
p^{r} & ($q$: the power mean).
\end{subnumcases}
\end{theorem}
\noindent\textbf{Proof.}
We can prove this theorem in the same way as Theorem 3.1 in \cite{nielsen2009sided}.
By the definition of the $g$-Bregman divergence, we have 
\begin{align}
\label{sumofdivergence}
\sum_{\nu=1}^N \alpha_\nu B_F^g(p_\nu,q)=\sum_{\nu=1}^N \alpha_\nu \biggl(F(g(p_\nu)) - F(g(q)) - \langle g(p_\nu) - g(q), \nabla F(g(q))\rangle\biggr) 
\end{align}
Let $g(p')\eqdef \sum_{\nu=1}^N \alpha_\nu g(p_\nu)$.
The RHS of  (\ref{sumofdivergence}) yields
\begin{align}
\biggl(\sum_{\nu=1}^N \alpha_\nu F(g(p_\nu))-F(g(p'))\biggr) + \biggl(F(g(p'))-F(g(q))\biggr)-\sum_{\nu=1}^N\alpha_\nu \langle g(p_\nu) - g(q), \nabla F(g(q))\rangle \\ \nonumber
=\biggl(\sum_{\nu=1}^N \alpha_\nu F(g(p_\nu))-F(g(p'))\biggr) + \biggl(F(g(p'))-F(g(q)) -\langle \sum_{\nu=1}^N\alpha_\nu g(p_\nu) - g(q), \nabla F(g(q))\rangle\biggr) \\ \nonumber
=\biggl(\sum_{\nu=1}^N \alpha_\nu F(g(p_\nu))-F(g(p'))\biggr) + B_F^g(p',q)\geq \sum_{\nu=1}^N \alpha_\nu F(g(p_\nu))-F(g(p')) \\ \nonumber
=\sum_{\nu=1}^N \alpha_\nu F(g(p_\nu))-F\biggl(\sum_{\nu=1}^N \alpha_\nu g(p_\nu) \biggr)=J_{F,\boldsymbol{\alpha}}^g(p_1,p_2,\cdots, p_N) .
\end{align}
Because $B_F^g(p',q)$ equal to zero if and only if $q=p'=g^{-1}\biggl(\sum_{\nu=1}^N \alpha_\nu g(p_\nu) \biggr)$, we prove the theorem.

\begin{prop}
\label{exchange_variable}
For the $g$-Bregman divergence, the following equation holds.
\begin{align}
B_F^g(q,p)=B_{F^\ast}^{\hat{g}}(p,q),
\end{align}
where $F^\ast(x^\ast)\eqdef \sup_x\{\langle x^\ast, x\rangle - F(x)\}$ denotes the Legendre convex conjugate and $\hat{g}\eqdef \nabla F \circ g$.
\end{prop}
\noindent\textbf{Proof.}
For the Bregman divergence, the equation $B_F(q',p')=B_F^{\ast}(\nabla F(p'),\nabla F(q'))$ holds\cite{nielsen2009sided,nielsen2011burbea}.
By putting $p'=g(p)$ and $q'=g(q)$, the result follows.

\begin{corollary}
\label{left_centroid}
Let $\Omega$ be a $g$-convex set.
Let $p_\nu(\nu=1,2,3\cdots,N)$ and $q$ are points in $\Omega$.
Let $\alpha_\nu\geq 0(\nu=1,2,3\cdots,N)$ be parameters which satisfy $\sum_{\nu=1}^N\alpha_\nu=1$.

Then, the following inequality holds.
\begin{align}
\sum_{\nu=1}^N \alpha_\nu B_F^g(q,p_\nu)\geq J_{F^\ast,\boldsymbol{\alpha}}^{\hat{g}}(p_1,p_2,\cdots, p_N)
\end{align}
Equality holds if and only if
\begin{align}
q=\hat{g}^{-1}\biggl(\sum_{\nu=1}^N \alpha_\nu\hat{g}(p_\nu) \biggr)
\end{align}
\end{corollary}
\noindent\textbf{Proof.}
From Proposition \ref{exchange_variable}, we obtain
\begin{align}
\sum_{\nu=1}^N \alpha_\nu B_F^g(q,p_\nu)=\sum_{\nu=1}^N \alpha_\nu B_{F^\ast}^{\hat{g}}(p_\nu,q)
\end{align}
Applying Theorem \ref{centroid_theorem} to this equation, the result follows.

From Theorem \ref{centroid_theorem} and Corollary \ref{left_centroid}, we conclude the point which minimize the weighted average of the $g$-Bregman divergence is unique.

\subsection{Relation between the $g$-Bregman and skew $g$-Jensen divergence}
In this subsection, we derive inequality between the $g$-Bregman and the skew $g$-Jensen divergence.
The following equations also hold for the $g$-Bregman and the skew $g$-Jensen divergence.

\begin{align}
B_F^g(q,p)&=\lim_{\alpha\downarrow 0}sJ_{F,\alpha}^g(p,q)\\
B_F^g(p,q)&=\lim_{\alpha\uparrow 1}sJ_{F,\alpha}^g(p,q)
\end{align}
It is easily proved by putting $p'=g(p)$ and $q'=g(q)$ in (\ref{jensen_lim1}) and (\ref{jensen_lim2}).

Then, we show some lemmas.
\begin{lemma}
\label{relation_Jensen_Bregman}
Let $\Omega$ be a $g$-convex set and let $p,q$ be points $p,q\in \Omega$. 
For a parameter $\alpha\in (0,1)$ and a point $r\in\Omega$ which satisfies $g(r)=(1-\alpha)g(p)+\alpha g(q)$, the $g$-Bregman divergence can be expressed as follows.
\begin{equation}
\label{g_Jensen_Bregman_eq}
sJ_{F,\alpha}^g(p ,q)\eqdef  \frac{1}{\alpha(1-\alpha)}\biggl((1-\alpha)B_{F}^g(p,r)+\alpha B_{F}^g(q,r)\biggr)
\end{equation}
\end{lemma}
It is easily proved in the same way as equation (\ref{Jensen_Bregman_eq}) by putting $p'=g(p)$ and $q'=g(q)$.

\begin{lemma}
\label{division}
Let $\Omega$ be a $g$-convex set and let $p,q$ be points $p,q\in \Omega$. 
For a parameter  $\alpha\in (0,1)$ and a point $r\in\Omega$ which satisfies $g(r)=(1-\alpha)g(p)+\alpha g(q)$, the following equations hold.
\begin{eqnarray}
\label{eq_division}
B_{F}^g(p,q)=B_{F}^g(p,r)+B_{F}^g(r,q)+\frac{\alpha}{1-\alpha}B_{F,\mathrm{sym}}^g(r,q)\\
\label{eq_inv_division}
B_{F}^g(q,p)=B_{F}^g(r,p)+B_{F}^g(q,r)+\frac{1-\alpha}{\alpha}B_{F,\mathrm{sym}}^g(r,p)
\end{eqnarray}
\end{lemma}
\noindent\textbf{Proof.}
We first prove (\ref{eq_division}).
By the generalized law of cosines (\ref{triangular relation}), we have
\begin{equation}
\label{line}
B_{F}^g(p,q)=B_{F}^g(p,r)+B_{F}^g(r,q)+\langle g(r)-g(p), \nabla F(g(q))- \nabla F(g(r)).
\end{equation}
From the assumption, the equation
\begin{equation}
\label{eq_ratio}
g(r)-g(p)=\frac{\alpha}{1-\alpha}(g(q)-g(r))
\end{equation}
holds.
By substituting (\ref{eq_ratio}) to (\ref{line}) and using (\ref{def_sym_breg}), the result follows. By exchanging $p$ and $q$ in (\ref{line}), we can prove (\ref{eq_inv_division}) in the same way.
\\

\begin{lemma}
\label{geometry_division}
Let $\Omega$ be a $g$-convex set and let $p,q$ be points $p,q\in \Omega$. 
For a parameter  $\alpha\in (0,1)$ and a point $r\in\Omega$ which satisfies $g(r)=(1-\alpha)g(p)+\alpha g(q)$, the following equation holds.
\begin{equation}
\label{eq_geometry_division}
B_{F,\mathrm{sym}}^g(p,q)=\frac{1}{\alpha}B_{F,\mathrm{sym}}^g(p,r)+\frac{1}{1-\alpha}B_{F,\mathrm{sym}}^g(r,q)
\end{equation}
\end{lemma}
\noindent\textbf{Proof.}
Taking the sum of (\ref{eq_division}) and (\ref{eq_inv_division}), the result follows.\\

\begin{theorem}(Bregman-Jensen inequality)\\
\label{Th_JB_inequality}
Let $\Omega$ be a $g$-convex set and let $p,q$ be points in $\Omega$.
For a parameter $\alpha\in(0,1)$, the skew $g$-Jensen divergence and the symmetric $g$-Bregman divergence, the following inequality holds.
\begin{equation}
\label{JB_inequality}
B_{F,\mathrm{sym}}^g(p,q)\geq sJ_{F,\alpha}^g(p ,q)
\end{equation}
\end{theorem}
\noindent\textbf{Proof.}\\
Let $r\in \Omega$ be a point which satisfies $g(r)=(1-\alpha)g(p)+\alpha g(q)$.
From (\ref{g_Jensen_Bregman_eq}) and using $ B_{F}^g(p,q)\leq B_{F,\mathrm{sym}}^g(p,q)$, we have 
\begin{equation}
sJ_{F,\alpha}^g(p ,q)= \frac{1}{\alpha}B_{F}^g(p,r)+\frac{1}{1-\alpha} B_{F}^g(q,r)\leq\frac{1}{\alpha}B_{F,\mathrm{sym}}^g(p,r)+\frac{1}{1-\alpha} B_{F,\mathrm{sym}}^g(q,r).
\end{equation}
Using Lemma \ref{geometry_division}, we have
\begin{align}
sJ_{F,\alpha}^g(p ,q)\leq \frac{1}{\alpha}B_{F,\mathrm{sym}}^g(p,r)+\frac{1}{1-\alpha}B_{F,\mathrm{sym}}^g(r,q)
=B_{F,\mathrm{sym}}^g(p,q)
\end{align}

\begin{corollary}(Bregman-Jensen inequality for parallelograms)\\
Let $\Omega$ be a $g$-convex set and let $p,q,r,s\in\Omega$ be points which satisfy $g(p)+g(r)=g(q)+g(s)$ (The points $p,q,r,s$ are vertices of ''parallelogram'').

Let $\boldsymbol{\alpha}=(\frac{1}{4},\frac{1}{4},\frac{1}{4},\frac{1}{4})$.

Then, the following inequality holds.
\begin{align}
\frac{1}{8}\biggl(B_{F,\mathrm{sym}}^g(p,q)+B_{F,\mathrm{sym}}^g(q,r)+B_{F,\mathrm{sym}}^g(r,s)+B_{F,\mathrm{sym}}^g(s,p)\biggr)\geq J_{F,\boldsymbol{\alpha}}^g(p,q,r,s)
\end{align}
\end{corollary}
\noindent\textbf{Proof.}
Let $c\in\Omega$ be a point which satisfies $g(c)=\frac{1}{2}\bigl(g(p)+g(r)\bigr)=\frac{1}{2}\bigl(g(q)+g(s)\bigr)$.
From Theorem \ref{Th_JB_inequality}, we have
\begin{align}
\label{BJ_1}
B_{F,\mathrm{sym}}^g(p,r)\geq sJ_{F,\frac{1}{2}}^g(p ,r)&=4\biggl(\frac{1}{2}F(g(p))+\frac{1}{2}F(g(r))-F(g(c))\biggr)\\ 
\label{BJ_2}
B_{F,\mathrm{sym}}^g(q,s)\geq sJ_{F,\frac{1}{2}}^g(q ,s)&=4\biggl(\frac{1}{2}F(g(q))+\frac{1}{2}F(g(s))-F(g(c))\biggr)
\end{align}
Taking the sum of (\ref{BJ_1}) and (\ref{BJ_2}) and applying Theorem \ref{parallelogram}, we have 
\begin{align}
B_{F,\mathrm{sym}}^g(p,q)+B_{F,\mathrm{sym}}^g(q,r)+B_{F,\mathrm{sym}}^g(r,s)+B_{F,\mathrm{sym}}^g(s,p)\\ \nonumber
\geq 2\biggl(F(g(p))+F(g(q))+F(g(r))+F(g(s))\biggr)-8F(g(c))
\end{align}
By the assumption, we get $g(c)=\frac{1}{4}\bigl(g(p)+g(q)+g(r)+g(s)\bigr)$.
By the definition of the multivariate skew $g$-Jensen divergence (Definition \ref{definition of multivariate}), the result follows.
\section{Examples}
We show some examples of the $g$-divergences.

We denote  by ''GBD'' the $g$-Bregman divergence, by ''SGBD'' the symmetric $g$-Bregman divergence, ''SGJD'' by the skew $g$-Jensen divergence and by ''BJ-inequality'' the Bregman-Jensen inequality (\ref{JB_inequality}).

A parameter $\alpha$ is in $(0,1)$ except for example 3) and we use a parameter $\beta\in(0,1)$ instead of $\alpha$ in example 3). The variables $p_i$ and $q_i$ are non-negative real numbers for all $i$. 

From example 1) to 6), we show the case that GBD belong to a class of $f$-divergence and example 7) is the case of information geometry\cite{amari2016information,amari2010information,amari2007methods} and statistical mechanics.\\

\textbf{1) KL-divergence, Jeffreys divergence, JS-divergence} \\
\begin{itemize}
\item Generating functions: $F(p)=\sum_i p_i\ln p_i$, $g(p)=p$\\
\item GBD: generalized KL-divergence $KL(p\|q)\eqdef \sum_i (-p_i+q_i+p_i\ln(\frac{p_i}{q_i}))$\\
\item SGBD: Jeffreys divergence (J-divergence)\cite{jeffreys1946invariant} $J(p,q)\eqdef\sum_i(p_i-q_i)\ln(\frac{p_i}{q_i})$ \\
\item SGJD: scaled skew JS-divergence \\$\frac{1}{\alpha(1-\alpha)}JS_{\alpha}(p\|q)\eqdef  \frac{1}{\alpha(1-\alpha)}\biggl((1-\alpha)KL(p\| (1-\alpha)p+\alpha q)+\alpha KL(q\| (1-\alpha)p+\alpha q)\biggr)$  \\
\item BJ-inequality: generalized Lin's inequality $\alpha(1-\alpha)J(p,q)\geq JS_{\alpha}(p,q)$. When $\alpha$ is equall to $\frac{1}{2}$, this inequality is consistent with Lin's inequality\cite{lin1991divergence}.\\ 
\end{itemize}

 \textbf{2) KL-divergence, Jeffreys divergence, $\alpha$-divergence} \\
\begin{itemize}
\item Generating functions: $F(p)=\sum_i \exp(p_i)$, $g(p)=\ln p$\\
\item GBD: generalized reverse KL-divergence $KL(q\|p)\eqdef \sum_i (p_i-q_i+q_i\ln(\frac{q_i}{p_i}))$\\
\item SGBD: J-divergence  $J(p,q)\eqdef\sum_i(p_i-q_i)\ln(\frac{p_i}{q_i})$\\
\item SGJD: generalized $\alpha$-divergence $D_{1-\alpha}(p\| q)$, \\
where $D_{\alpha}(p\| q)\eqdef \sum_i \frac{1}{\alpha(\alpha-1)}(p_i^\alpha q_i^{1-\alpha} -\alpha p_i -(1-\alpha)q_i)$\\
\item BJ-inequality: $J(p,q)\geq D_{\alpha}(p\| q)$\\
\end{itemize}

\textbf{3) $\alpha$-divergence} \\
\begin{itemize}
\item Generating functions: $F(p)=\frac{1}{1-\alpha}\sum_i p_i^{\frac{1}{\alpha}}$, $g(p)=p^\alpha$ for $\alpha\in \mathbb{R}\setminus\{0,1\}$\\
\item GBD:  generalized $\alpha$-divergence $D_{\alpha}(p\|q)$\\
\item SGBD: symmetric  $\alpha$-divergence $D_{\alpha,\mathrm{sym}}(p,q)\eqdef D_{\alpha}(p\|q)+D_{\alpha}(q\|p)$\\
\item SGJD: $\frac{1}{\beta(1-\beta)(1-\alpha)}\sum_i\biggl((1-\beta)p_i +\beta q_i -{\bigl((1-\beta)p_i^\alpha+\beta q_i^\alpha\bigr)}^\frac{1}{\alpha}\biggr)$\\
\item BJ-inequality: $D_{\alpha,\mathrm{sym}}(p,q)\geq\frac{1}{\beta(1-\beta)(1-\alpha)}\sum_i\biggl((1-\beta)p_i +\beta q_i -{\bigl((1-\beta)p_i^\alpha+\beta q_i^\alpha\bigr)}^\frac{1}{\alpha}\biggr)$\\
\end{itemize}
By using the equations  $\frac{1}{2}D_{(\frac{1}{2})}(p\|q)=H^2(p,q)$, $2D_{(2)}(p\|q)=D_{\chi^2,P}(p\|q)$ and $2D_{(-1)}(p\|q)=D_{\chi^2,N}(p\|q)$ \cite{cichocki2010families}, we obtain example 4) to 6).\\ 

 \textbf{4) Hellinger distance} \\
\begin{itemize}
\item Generating functions: $F(p)=\sum_i p_i^2$, $g(p)=\sqrt{p}$\\
\item GBD: Hellinger distance $H^2(p,q)\eqdef\sum_i(\sqrt{p_i}-\sqrt{q_i})^2 $\\
\item SGBD: Hellinger distance $2H^2(p,q)$\\
\item SGJD: Hellinger distance $H^2(p,q)$ \\
\item BJ-inequality: trivial\\
\end{itemize}

\textbf{5) Pearson $\chi$-square divergence} \\
\begin{itemize}
\item Generating functions: $F(p)=-2\sum_i \sqrt{p_i}$, $g(p)=p^2$\\
\item GBD:  Pearson $\chi$-square divergence $D_{\chi^2,P}(q\|p)\eqdef\sum_i \frac{(p_i-q_i)^2}{q_i}$\\
\item SGBD: symmetric $\chi$-square divergence $D_{\chi^2,\mathrm{sym}}(p,q)\eqdef D_{\chi^2}(p\|q)+D_{\chi^2}(q\|p)$\\
\item SGJD: $\frac{2}{\alpha(1-\alpha)}\sum_i(-(1-\alpha)p_i - \alpha q_i + \sqrt{(1-\alpha)p_i^2+\alpha q_i^2})$\\
\item BJ-inequality: $D_{\chi^2,\mathrm{sym}}(p,q)\geq \frac{2}{\alpha(1-\alpha)}\sum_i(-(1-\alpha)p_i - \alpha q_i + \sqrt{(1-\alpha)p_i^2+\alpha q_i^2})$\\
\end{itemize}

\textbf{6) Neyman $\chi$-square divergence} \\
\begin{itemize}
\item Generating functions: $F(p)=\sum_i {p_i}^{-1}$, $g(p)=p^{-1}$\\
\item GBD:  Neyman $\chi$-square divergence $D_{\chi^2,N}(q\|p)\eqdef\sum_i \frac{(p_i-q_i)^2}{p_i}$\\
\item SGBD: symmetric $\chi$-square divergence $D_{\chi^2,\mathrm{sym}}(p,q)$\\
\item SGJD: $\frac{1}{\alpha(1-\alpha)}\sum_i((1-\alpha)p_i + \alpha q_i - \frac{p_iq_i}{(1-\alpha)q_i+\alpha p_i})$\\
\item BJ-inequality: $D_{\chi^2,\mathrm{sym}}(p,q)\geq \frac{1}{\alpha(1-\alpha)}\sum_i((1-\alpha)p_i + \alpha q_i - \frac{p_iq_i}{(1-\alpha)q_i+\alpha p_i})$\\
\end{itemize}

\textbf{7) Information geometry and statistical mechanics}\\ 
In the dually flat space $M$, there exist two affine coordinates $\theta^i$, $\eta_i$ and two convex functions called potential $\psi(\theta)$, $\phi(\eta)$.
For the points $P,Q\in M$, we put $p_i=\theta^i(P)$, $q_i=\theta^i(Q)$ or  $p_i=\eta_i(P)$, $q_i=\eta_i(Q)$.\\

\begin{itemize}
\item Generating functions: $F(p)=\psi(\theta)$ or $F(p)=\phi(\eta)$, $g(p)=p$\\
\item GBD: canonical divergence $D(P\|Q)\eqdef \phi(\eta(Q)) + \psi(\theta(P)) - \sum_i \eta_i(Q)\theta^i(P)$\\
\item SGBD: $D_A(P,Q)\eqdef \sum_i(\eta_i(Q)-\eta_i(P))(\theta^i(Q)-\theta^i(P))$\\
\item SGJD: $sJ_{\psi,\alpha}(\theta(P),\theta(Q))$ or  $sJ_{\phi,\alpha}(\eta(P),\eta(Q))$\\
\item BJ-inequality: $D_A(P,Q)\geq sJ_{\psi,\alpha}(\theta(P),\theta(Q))$ or $D_A(P,Q)\geq sJ_{\phi,\alpha}(\eta(P),\eta(Q))$\\
\end{itemize}
For the exponential or mixture families, GBD is equal to the KL-divergence and SGBD is equal to the J-divergence.
For the exponential families, SGJD $sJ_{\psi,\alpha}(\theta(P),\theta(Q))$ is equal to the skew Bhattacharyya distance\cite{nielsen2011burbea,bhattacharyya1943measure }, and for the mixture families , SGJD $sJ_{\phi,\alpha}(\eta(P),\eta(Q))$ is equal to the skew-JS divergence(see \cite{nishiyama2018divergence} in detail).

In the case of canonical ensemble in statistical mechanics, the probability distribution function belongs to exponential families. Because the manifold of the exponential family is a dually flat space, the same equations and inequality hold for 
\begin{align}
\theta&= -\beta\equiv -\frac{1}{k_BT} \\
\eta&=U \\
\psi&= -\beta F\\
\phi&=-\frac{1}{k_B}S,
\end{align}
where $k_B$ is the Boltzmann constant, $T$ is temperature, $U$ is the internal energy, $S$ is the entropy and $F$ is the Helmholtz free energy.

Defining a parameter $\alpha$ as a index of the state which satisfy $\frac{1}{T_\alpha}=(1-\alpha)\frac{1}{T_0}+\alpha\frac{1}{T_1}$ or  $U_\alpha=(1-\alpha)U_0+\alpha U_1$, BJ-inequality can be written as follows.
\begin{align}
\alpha(1-\alpha)\frac{(T_1-T_0)(U_1-U_0)}{T_0T_1}\geq \frac{F_\alpha}{T_\alpha}-(1-\alpha)\frac{F_0}{T_0}-\alpha\frac{F_1}{T_1}\\
\alpha(1-\alpha)\frac{(T_1-T_0)(U_1-U_0)}{T_0T_1} \geq S_\alpha-(1-\alpha)S_0-\alpha S_1
\end{align}

\section{Conclusion}
We have introduced the $g$-Bregman divergence and the skew $g$-Jensen divergence by extending the definitions of the Bregman and the skew Jensen divergence respectively. 

First, we have shown the geometrical properties of the $g$-Bregman divergence and the symmetric $g$-Bregman divergence.

Then, we have derived an inequality between the symmetric $g$-Bregman divergence and the skew $g$-Jensen divergence.

Finally, we have shown they include divergences which belong to a class of $f$-divergence .

If the relationship between these three classes of divergence are studied in more detail, it is expected to proceed applications to various fields including machine learning.

\bibliography{reference_v6}
\end{document}